\documentclass{elsart}
\usepackage{graphics}
\usepackage{graphicx}
\usepackage{epsfig}
\usepackage{amssymb,amsmath}
\usepackage[colorlinks]{hyperref}

\begin{document}
\begin{frontmatter}
\title{Some Ultraspheroidal Monogenic Clifford Gegenbauer Jacobi Polynomials and Associated Wavelets}
\author{Sabrine Arfaoui}
\address{Department of Informatics, Higher Institute of Applied Sciences and Technology of Mateur, Street of Tabarka, 7030 Mateur, Tunisia.\\
and\\
Research Unit of Algebra, Number Theory and Nonlinear Analysis UR11ES50, Faculty of Sciences, Monastir 5000, Tunisia.}
\ead{arfaoui.sabrine@issatm.rnu.tn}
\author{Anouar Ben Mabrouk}
\address{Higher Institute of Applied Mathematics and Informatics, University of Kairouan, Street of Assad Ibn Alfourat, Kairouan 3100, Tunisia.\\
and\\
Research Unit of Algebra, Number Theory and Nonlinear Analysis UR11ES50, Faculty of Sciences, Monastir 5000, Tunisia.}
\ead{anouar.benmabrouk@fsm.rnu.tn}
\begin{abstract}
In the present paper, new classes of wavelet functions are presented in the framework of Clifford analysis. Firstly, some classes of orthogonal polynomials are provided based on 2-parameters weight functions. Such classes englobe the well known ones of Jacobi and Gegenbauer polynomials when relaxing one of the parameters. The discovered polynomial sets are next applied to introduce new wavelet functions. Reconstruction formula as well as Fourier-Plancherel rules have been proved.
\end{abstract}
\end{frontmatter}
\section{Introduction}
Spheroidal functions or precisely spheroidal wave functions are issued from the wave equation
$$
\nabla^2w+k^2w=0.
$$
By considering solutions of separated variables in an elliptic cylinder coordinates system, or the prolate or oblate spheroids, such solutions satisfy a second order ODE of the form
$$
(1-t^2)w''+2\alpha tw'+(\beta-\gamma^2t^2)w=0.
$$
for both radial and angular functions. In fact, prolate and oblate spheroidal coordinate systems are results of rotating the two-dimensional elliptic coordinate system, consisting of confocal ellipses and hyperbolas, about the major and minor axes of the ellipses. See \cite{Abramowitzetal}, \cite{Lietal}, \cite{Morais}, \cite{Stratton}, \cite{Strattonetal} \cite{Osipovetal}. This last equation leads to special functions such as Bessel, Airy, ... and special polynomials such as Gegenbauer, Legendre, Chebyshev, .... This is a first idea behind the link between these functions and a first motivation of our work and its titling as spheroidal wavelets. Besides, spheroidal functions have been in the basis of modeling physical phenomena where the wave behaviour is pointed out such as radars, antennas, 3D-images, ... Recall also that Gegenbauer polynomials themselves are strongly related to spheroidal functions since their appearance and these are called ultraspheroidal polynomials. See \cite{Antoine-Murenzi-Vandergheynst}, \cite{DeSchepper}, \cite{Delanghe}, \cite{Lehar}, \cite{Lietal}, \cite{Michel}, \cite{Moussa}, \cite{Saillardetal}.

The use of wavelets in the analysis of functions is widespread especially in the last decades. Nowadays, wavelets are interesting and useful tools in many fields such as mathematics, quantum physics, electrical engineering, time/image processing, bio-signals, seismology, geology, .....
	
Wavelets have been created to meet a need in signal processing that is not well understood by Fourier theory. Classical Fourier analysis provides a global approach for signals as it replaces the analyzed function with a whole-space description (See (\ref{Fourier-transform-onRm}) later). Wavelet analysis in contrast decomposes the signal in both time and frequency and describes it locally and globally, as the need.
	
Wavelet analysis of a function $f$ in the space of able analyzed functions (generally $L_2$) starts by convoluting it with local copy of a wavelet mother function $\psi$ known as the analyzing wavelet relatively to 2-parameters; One real number parameter $a>0$ defines the dilation parameter or the scale and one space parameter $b$ in the same space as the function $f$ and $\psi$ domains defines the translation parameter or the position. Such copy is denoted usually by $\psi_{a,b}$ and is defined by
\begin{equation}\label{psiab}
\psi_{a,b}(x)=a^{-\frac{1}{2}}\psi(\displaystyle\frac{x-b}{a}).
\end{equation}
To be a good candidate as a wavelet mother, an admissibility assumption on the function $\psi$ is usually assumed. It states that
\begin{equation}\label{admissibility-conditionofpsi}
\mathcal{A}_\psi=\displaystyle\int_{-\infty}^{+\infty}\displaystyle\frac{|\widehat{\psi}(u)|^2}{|u|}du<+\infty,
\end{equation}
where $\widehat{\psi}$ is the Fourier transform of $\psi$.
	
The convolution of the analyzed function $f$ with the copy $\psi_{a,b}$ defines the so-called wavelet transform of $f$ or exactly the Continuous Wavelet Transform (CWT) expressed by
\begin{equation}\label{waveletcoefficientcab(f)}
C_{a,b}(f)=<f,\psi_{a,b}>=\,\displaystyle\int_{-\infty}^{+\infty}f(x)
\overline{\psi_{a,b}(x)}dx.
\end{equation}
Whenever the admissibility condition is fulfilled, the analyzed function $f$ may be reconstructed in an $L_2$ sense as
\begin{equation}
f(x)=\displaystyle\frac{1}{\mathcal{A}_\psi}\displaystyle\int_{\mathbb{R}}\displaystyle\int_{0}^{+\infty} C_{a,b}(f)\psi_{a,b}(x)\displaystyle\frac{da}{a^2}db,
\end{equation}
where the equality has to be understood in the $L_2$-sense (See \cite{Jaffard1}, \cite{Jaffard2}). This equality will be proved later in the present context of Clifford Gegenbauer-Jacobi type wavelets.
	
Usually analyzing wavelets are related also to moments. The regularity of the analyzing wavelet $\psi$ is related to a number of vanishing moments that should be satisfied
\begin{equation}\label{vanishingmomentsofpsi}
\displaystyle\int_{-\infty}^{+\infty} x^n \psi(x) dx=0,\quad n=0,1,\dots,N.
\end{equation}
Such a condition helps to analyze functions of some fixed regularity. In wavelet theory, the first result relating regularity to wavelet transforms is due to Jaffard (See \cite{Holschneider-Tchamitchan}, \cite{Jaffard1}, \cite{Jaffard2}) and is stated as follows.
\begin{prop}\label{PropJaddard} Let $\psi$ be a $C^{r}(\mathbb{R}^m)$ function with all moments of order less than $r$ vanishing and all derivatives of order less than $r$ well localized.
\begin{itemize}
\item $f\in\mathcal{C}^{\alpha}(\mathbb{R}^m)$ if and only if $|C_{a,b}(f)|\leq\,Ca^{\alpha}$ for all $b$ and $0<a<<1$.
\item If $f\in\mathcal{C}^{\alpha}(x_0)$, then for $0<a<<1$ and $|b-x_0|\leq1/2$,
\begin{equation}\label{ch3:coeff}
|C_{a,b}(f)|\leq\,Ca^{\alpha}\left(1+\frac{|b-x_0|}{a}\right)^{\alpha}.
\end{equation}
\item If (\ref{ch3:coeff}) holds and if $f\in\mathcal{C}^\varepsilon(\mathbb{R}^m)$ for an $\varepsilon>0$, then there exists a polynomial $P$ such that, if $|x-x_0|\leq1/2$,
\begin{equation}\label{ch3:eq:1.11}
|f(x)-P(x-x_0)|\leq\,C|x-x_0|^{\alpha}\log\left(\frac{2}{|x-x_0|}\right).
\end{equation}
\end{itemize}
\end{prop}
More about regularity, admissibility, vanishing moments and wavelet properties may be found in \cite{Craddocketal}, \cite{Delanghe}, \cite{Kilbas}, \cite{Mitrea}, \cite{Pena}, \cite{Vieira}.

It holds that wavelet theory on the real line and generally on Euclidian spaces has been extended in some cases of Clifford analysis. The classical wavelet theory can be constructed in the framework of Clifford analysis. Clifford analysis deals with so-called monogenic functions which are described as solutions of the Dirac operator and/or direct higher dimensional generalizations of holomorphic functions in the complex plane. Clifford wavelets and the possibility to construct orthogonal wavelet bases and consequently multiresolution analyses associated has been the object of several works, but remain to be a fascinating subject of researches. In \cite{Askarietal} a multiresolution analysis in the context of Clifford analysis has been provided. Clifford scaling functions, Clifford wavelets as well as related wavelet filters has been developed.and proved to be applicable in quantum mechanics. In \cite{Kumar1} and \cite{Kumar2}, spheroidal wavelets leading to frames as well as multiresolution analysis have been developed. It was proved that spheroidal functions may induce good candidates characterized by localizations in both frequency and space and thus lead to good wavelets. More facts about Clifford wavelets and discussions on possible associated multiresolution analyses may be found in \cite{Brackx-Schepper-Sommen0}, \cite{Brackx-Schepper-Sommen1}, \cite{Brackx-Schepper-Sommen2}, \cite{Brackx-Schepper-Sommen3}, \cite{Brackx-Schepper-Sommen4}, \cite{Brackx-Schepper-Sommen5}, \cite{Hitzeretal}.
	
Let $\Omega$ be an open subset of $\mathbb{R}^m$ or $\mathbb{R}^{m+1}$ and $f:\Omega\rightarrow\mathbb{A}$, where $\mathbb{A}$ is the real Clifford algebra $\mathbb{R}_{m}$ (or $\mathbb{C}_{m}$). $f$ may be written in the form
\begin{equation}\label{fincliffordalbebra}
f=\displaystyle\sum_{A}f_{A}e_{A}
\end{equation}
where the functions $f_A$ are $\mathbb{R}$ (or $\mathbb{C}$)-valued and $(e_A)_A$ is a suitable basis of $\mathbb{A}$.
	
Despite the fact that Clifford analysis generalizes the most important features of classical complex analysis, monogenic functions do not enjoy all properties of holomorphic functions of one complex variable. For instance, due to the non-commutativity of the Clifford algebras, the product of two monogenic functions is in general not monogenic. It is therefore natural to look for specific techniques to construct monogenic functions. See \cite{Brackx-Schepper-Sommen0}, \cite{Delanghe}, \cite{Pena}.
	
In the literature, there are several techniques available to generate monogenic functions such as the Cauchy-Kowalevski extension (CK-extension) which consists in finding a monogenic extension $g^*$ of an analytic function $g$ defined on a given subset in $\mathbb{R}^{m+1}$ of positive codimension. For analytic functions $g$ on the plane $\{(x_0, \underline{x}) \in\mathbb{R}^{m+1}, \quad x_0 = 0\}$ the problem may be stated as follows: \textit{Find $g^*\in\mathbb{A}$ such that
\begin{equation}\label{cauchy-kowalevski-extension}
\partial_{x_0}g^*=-\partial_{\underline{x}} g^*\quad in \quad \mathbb{R}^{m+1}\quad\hbox{and}\quad
g^*(0,\underline{x})=g(\underline{x}).
\end{equation}}
A formal solution is
\begin{equation}\label{ck-formal-solution}
g^*(x_0,\underline{x})=\exp(-x_0\partial_{\underline{x}}) g(\underline{x})=\displaystyle\sum_{k=0}^{\infty}\displaystyle\frac{(-x_0)^k}{k!}
\partial_{\underline{x}}^kg(\underline{x}).
\end{equation}
It may be proved that (\ref{ck-formal-solution}) is a monogenic extension of the function $g$ in $\mathbb{R}^{m+1}$. Moreover, by the uniqueness theorem for monogenic functions this extension is also unique. See \cite{Brackx-Schepper-Sommen0}, \cite{Delanghe}, \cite{Pena}, \cite{Vieira}, \cite{Winkler} and the references therein.
		
The organization of this paper is as follows: In section 2, a brief overview of some properties of the Clifford and Fourier analysis has been conducted. Section 3 is devoted to a review of the class of Gegenbauer-Jacobi polynomials in the framework of Clifford analysis. In section 4, some new classes of polynomials generalizing those of section 3 are developed by adapting 2-parameters weights and thus applied to introduce some new wavelets. Section 5 is devoted to the link and discussions about the present case and Legendre and Tchebyshev polynomials as well as the role of the parameters $\alpha$ and $\beta$ in the Clifford weight function applied here. We concluded afterward.
\section{Clifford analysis revisited}
In this section we revisit some basic concepts that will be used later. Let $f$ be in $L^1(\mathbb{R}^m)$. Its Fourier transform denoted usually $\widehat{f}$ or $\mathcal{F}(f)$ is given by
\begin{equation}\label{Fourier-transform-onRm}
\widehat{f}(\eta)=\mathcal{F}(f)(\eta)=\displaystyle\frac{1}{(2\pi)^{\frac{m}{2}}}\displaystyle\int_{\mathbb{R}^m}\exp(-ix.\eta)f(x)dx,
\end{equation}
where $dx$ is the Lebesgues measure on $\mathbb{R}^m$ and $x.\eta$ is the standard inner product of $x$ and $\eta$ in $\mathbb{R}^m$.
		
Clifford analysis appeared as a generalization of the complex analysis and Hamiltonians. It extended complex calculus to some type of finite-dimensional associative algebra known as Clifford algebra endowed with suitable operations as well as inner products and norms. It is now applied widely in a variety of fields including geometry and theoretical physics. See \cite{Brackx-Schepper-Sommen0}, \cite{Delanghe}, \cite{Hitzeretal}, \cite{Lehar}, \cite{Lietal}, \cite{McIntoshetal1}, \cite{Pena}, \cite{Saillardetal}, \cite{Son}, \cite{Vieira} and the references therein.
		
Clifford analysis offers a functional theory extending the one of holomorphic functions of one complex variable. Starting from the real space $\mathbb{R}^m,\;(m\geq2)$ (or $\mathbb{C}^m$) endowed with an orthonormal basis $(e_1,\dots, e_m)$, the Clifford algebra $\mathbb{R}_m$ (or $\mathbb{C}_m$) starts by introducing a suitable interior product. Let
$$
e_j^2=-1,\quad j=1,\dots,m,
$$
$$
e_je_k+e_ke_j=0,\quad j\neq k,\quad j,k=1,\dots,m.
$$
It is straightforward that this is a non-commutative multiplication. Two anti-involutions on the Clifford algebra are important. The conjugation is defined as the anti-involution for which
$$
\overline{e_j}=-e_j,\quad j=1,\dots, m.
$$
The inversion is defined as the anti-involution for which
$$
e_j^{+}=e_j,\quad j=1,\dots,m.
$$
This yields a basis of the Clifford algebra ($e_A:A\subset\{1,\dots,m\}$)
where $e_{\emptyset}=1$
is the identity element. As these rules are defined, the Euclidian space $\mathbb{R}^m$ is then embedded in the Clifford algebras $\mathbb{R}_m$ (or $\mathbb{C}_m$) by identifying the vector $x=(x_1,\dots,x_m)$ with the vector $\underline{x}$ given by
$$
\underline{x}=\displaystyle\sum_{j=1}^{m}e_jx_j.
$$
The product of two vectors is given by
$$
\underline{x}\,\underline{y}=\underline{x}.\underline{y}+\underline{x}\wedge\underline{y}
$$
where
$$
\underline{x}.\underline{y}=-<\underline{x},\underline{y}>=-\displaystyle\sum_{j=1}^{m}x_j\,y_j
$$
and
$$
\underline{x}\wedge\underline{y}=\displaystyle\sum_{j=1}^{m}\displaystyle\sum_{k=j+1}^{m}e_j\,e_k(x_j\,y_k-x_ky_j).
$$
is the wedge product. In particular,
$$
\underline{x}^2=-<\underline{x},\underline{x}>=-|\underline{x}|^2.
$$
An $\mathbb{R}_m$ or $\mathbb{C}_m$-valued function $F(x_1,\dots,x_m)$, respectively $F(x_0, x_1,\dots,x_m)$ is called right monogenic in an open region of $\mathbb{R}^m$, respectively, or $\mathbb{R}^{m+1}$, if in that region
$$
F\partial_{\underline{x}}=0, \quad\mbox{respectively}\quad F(\partial_{x_0}+\partial_{\underline{x}})=0.
$$
Here $\partial_{\underline{x}}$ is the Dirac operator in $\mathbb{R}^m$ defined by
$$
\partial_{\underline{x}}=\displaystyle\sum_{j=1}^{m} e_j \partial_{x_j}
$$
and which splits the Laplacian in $\mathbb{R}^m$ as
$$
\Delta_m=-\partial_{\underline{x}}^2,
$$
whereas $\partial_{x_0}+\partial_{\underline{x}}$ is the Cauchy-Riemann operator in $\mathbb{R}^{m+1}$ for which
$$
\Delta_{m+1}=(\partial_{x_0}+\partial_{\underline{x}})(\partial_{x_0}+\overline{\partial_{\underline{x}}})
$$
Introducing spherical co-ordinates in $\mathbb{R}^m$ by
$$
\underline{x}=r\underline{\omega},\quad r=|\underline{x}|\in[0,+\infty[,\,\underline{\omega}\in S^{m-1},
$$
where $S^{m-1}$ is the unit sphere in $\mathbb{R}^m$, the Dirac operator takes the form
$$
\partial_{\underline{x}}=\underline{\omega}\left( \partial_r+\displaystyle\frac{1}{r} \Gamma_{\underline{\omega}}\right)
$$
where
$$
\Gamma_{\underline{\omega}}=-\displaystyle\sum_{i<j}e_ie_j(x_i\partial_{x_j}-x_j\partial_{x_i})
$$
is the so-called spherical Dirac operator which depends only on the angular co-ordinates.
		
As for the Euclidian case, Fourier analysis is extended to Clifford Fourier analysis \cite{Brackx-Schepper-Sommen0}, \cite{Brackx-Schepper-Sommen2}, \cite{Brackx-Schepper-Sommen3}, \cite{Craddocketal}. The idea behind the definition of the Clifford Fourier transform originates from the operator exponential representation of the classical Fourier transform by means of Hermite operators. Throughout this article the Clifford-Fourier transform of $f$ is given by
$$
\mathcal{F}(f(x))(y)=\displaystyle\int_{\mathbb{R}^m} e^{-i<\underline{x},\underline{y}>}\, f(\underline{x}) dV(\underline{x}),
$$
where $dV(\underline{x})$ is the Lebesgue measure on $\mathbb{R}^m$.

In the present work, we propose to apply such topics to output some generalizations of multidimensional Continuous Wavelet Transform in the context of Clifford analysis.
\section{Some old orthogonal polynomials revisited}
Firstly, we stress on the fact that the results presented in this section are not purely new. The same problem is already studied in \cite{Brackx-Schepper-Sommen1}. (See also \cite{Brackx-Schepper-Sommen4}, \cite{Brackx-Schepper-Sommen5}, \cite{DeSchepper} for similar results)
		
We propose to review the context of real Gegenbauer polynomials on $\mathbb{R}$ which are associated in this case to the real weight function $\omega(x)=(1+x^2)^\alpha$, $\alpha\in\mathbb{R}$, to the context of Clifford algebra-valued polynomials by considering the same weight function on the Clifford algebra $\mathbb{R}_m$. So, consider the Clifford algebra-valued weight function
$$
\omega_\alpha(\underline{x})=(1+|\underline{x}|^2)^\alpha,\, \alpha\in \mathbb{R}.
$$
The general Clifford-Gegenbauer polynomials, denoted by  $G_{\ell,m,\alpha}(\underline{x})$, are generated by the CK-extension $F^{*}(t,\underline{x})$ defined by
$$
F^{*}(t,\underline{x})=\displaystyle\sum_{\ell=0}^{\infty} \displaystyle\frac{t^\ell}{\ell!}G_{\ell,m,\alpha}(\underline{x})\,\omega_{\alpha-\ell}(\underline{x});\;\;t\in\mathbb{R},\;\;\underline{x}\in\mathbb{R}_m.
$$
As for the real case of orthogonal polynomials, we impose a left monogenic property on $F^{*}$ in $\mathbb{R}^{m+1}$ to obtain a recursive relation on the general Clifford-Gegenbauer polynomials $G_{\ell,m,\alpha}$. Hence, $F^{*}$ is monogenic means that
\begin{equation}\label{monogenic-property-GCGP}
(\partial_{t}+ \partial_{\underline{x}}) F^{*}(t,\underline{x})=0.
\end{equation}
The first part related to the time derivative is evaluated as
$$
\partial_tF^*(t,\underline{x})=\displaystyle\sum_{\ell=0}^{\infty}\displaystyle\frac{t^\ell}{\ell!}G_{\ell+1,m,\alpha}(\underline{x}) \,\omega_{\alpha-\ell-1}(\underline{x}).
$$
\begin{lem}\label{lemmaevenodd}
The Dirac operator of $\underline{x}^n$ is given by
\begin{equation}
\partial_{\underline{x}}(\underline{x}^n)=\gamma_{n,m} \underline{x}^{n-1}.
\end{equation}
where,
$$
\gamma_{n,m}=
\begin{cases}
-n \quad \hbox{if} \quad n\, \mbox{is \,even}.\\
-(m+n-1) \quad \mbox{if}\quad n \,\mbox{is\, odd}.
\end{cases}
$$
\end{lem}
\hskip-20pt Now, observing that
$$
\partial_{\underline{x}}(\underline{x})=-m,\quad\partial_{\underline{x}}(\underline{x}^2)=-2\underline{x}\quad\hbox{and}\quad
\partial_{\underline{x}}(|\underline{x}|^2)=2\underline{x},
$$
we get
$$
\partial_{\underline{x}}F^*(t,\underline{x})=\displaystyle\sum_{\ell=0}^{\infty}\displaystyle\frac{t^\ell}{\ell!}\left(\partial_{\underline{x}} G_{\ell,m,\alpha}(\underline{x})\omega_{\alpha-\ell}(\underline{x})+G_{\ell,m,\alpha}(\underline{x}) \partial_{\underline{x}}\omega_{\alpha-\ell}(\underline{x})\right).
$$
Observing again that
$$
\partial_{\underline{x}}\omega_{\alpha-\ell}(\underline{x})=2(\alpha-\ell)\underline{x}\,\,\omega_{\alpha-\ell-1}(\underline{x}),
$$
the monogenicity property (\ref{monogenic-property-GCGP}) leads to the recurrence relation
$$
\begin{array}{lll}
&&G_{\ell+1,m,\alpha}(\underline{x})\omega_{\alpha-\ell-1}(\underline{x})+\omega_{\alpha-\ell}(\underline{x})\partial_{\underline{x}}G_{\ell,m,\alpha}(\underline{x})\\
&&+2(\alpha-\ell) \omega_{\alpha-\ell-1}(\underline{x})\underline{x}G_{\ell,m,\alpha}(\underline{x})=0,
\end{array}
$$
or equivalently
\begin{equation}\label{rec}
G_{\ell+1,m,\alpha}(\underline{x})=-2(\alpha-\ell)\underline{x}G_{\ell,m,\alpha}(\underline{x})-(1+|\underline{x}|^2)
\partial_{\underline{x}}G_{\ell,m,\alpha}(\underline{x}).
\end{equation}
Starting from $G_{0,m,\alpha}(\underline{x})=1$, we obtain as examples
$$
G_{1,m,\alpha}(\underline{x})=-2\alpha \underline{x},
$$
$$
G_{2,m,\alpha}(\underline{x})=2\alpha[(2(\alpha-1)+m)\underline{x}^2-m],
$$
$$
G_{3,m,\alpha}(\underline{x})=[-4\alpha((2\alpha-1)+m)(\alpha-1)]\underline{x}^3+4\alpha(\alpha-1)(m+2)\underline{x}.
$$
The Clifford-Gegenbauer polynomials may be also introduced via the Rodrigues formula.
\begin{prop}
\begin{equation}\label{rod}
G_{\ell,m,\alpha}(\underline{x})=(-1)^{\ell}\,\,\omega_{\ell-\alpha}(\underline{x})\partial_{\underline{x}}^\ell(\,\omega_{\alpha}(\underline{x})).
\end{equation}
\end{prop}
\hskip-20pt\textbf{Proof.} We proceed by recurrence on $\ell$. For $\ell=1$, we have
$$
\partial_{\underline{x}}\,\omega_{\alpha}(\underline{x})=2\alpha\,\underline{x}\,\omega_{\alpha-1}(\underline{x})
=(-1)(-2\alpha\underline{x}) \,\omega_{\alpha-1}(\underline{x})=(-1)\,\omega_{\alpha-1}(\underline{x}) G_{1,m}^{\alpha}(\underline{x}).
$$
Which means that
$$
G_{1,m}^{\alpha}=(-1)\,\omega_{1-\alpha}(\underline{x}) \partial_{\underline{x}}\,\omega_{\alpha}(\underline{x}).
$$
For $\ell=2$, we get
$$
\begin{array}{lll}
\partial_{\underline{x}}^{(2)}\, \omega_{\alpha}  (\underline{x})&=&
2\alpha[ 2(\alpha-1)\underline{x}^2(1-\underline{x}^2)^{\alpha-2}-m(1-\underline{x}^2)^{\alpha-1} ]\\
&=& (-1)^2 \omega_{\alpha-2}[2\alpha[ 2(\alpha-1)+m]\underline{x}^2-m ]\\
&=& (-1)^2 \omega_{\alpha-2}(\underline{x}) \,G_{2,m,\alpha}(\underline{x}).
\end{array}
$$
Hence,
$$
G_{2,m}^{\alpha,\beta}=(-1)^2\omega_{2-\alpha}(\underline{x})\partial_{\underline{x}}^{(2)}\omega_{\alpha}(\underline{x}).
$$
So, assume that
$$
G_{\ell,m}^{\alpha,\beta}(\underline{x})=(-1)^\ell\omega_{\ell-\alpha}(\underline{x})\partial_{\underline{x}}^{(\ell)}\omega_{\alpha}(\underline{x}).
$$
Denote
$$
\Im(\underline{x})=-2(\alpha-\ell)\underline{x}(-1)^{\ell}\omega_{\ell-\alpha}(\underline{x}) \partial_{\underline{x}}^{(\ell)}\omega_{\alpha}(\underline{x}),
$$
and
$$
\Re(\underline{x})=(1+|\underline{x}|^2)(-1)^\ell2(\ell-\alpha)\underline{x} \omega_{\ell-\alpha-1}(\underline{x})\,\partial_{\underline{x}}^{\ell}\omega_{\alpha}(\underline{x}).
$$
From (\ref{rec}) and (\ref{rod}) we obtain
$$
\begin{array}{lll}
G_{\ell+1,m}^{\alpha,\beta}(\underline{x})&=&-2(\alpha-\ell)\underline{x}(-1)^{\ell}\omega_{\ell-\alpha}(\underline{x}) \partial_{\underline{x}}^{(\ell)}\omega_{\alpha}(\underline{x})\\
&\quad& -(1+|\underline{x}|^2)\partial_{\underline{x}}[(-1)^{\ell}\omega_{\ell-\alpha}(\underline{x}) \partial _{\underline{x}}^{(\ell)} \omega_{\alpha}(\underline{x})]\\
&=& \Im(\underline{x})- \Re(\underline{x})
-(1+|\underline{x}|^2)(-1)^{\ell}\omega_{\ell-\alpha}(\underline{x}) \partial_{\underline{x}}^{(\ell+1)} \omega_{\alpha}(\underline{x}).
\end{array}
$$
Simple calculus yield that
$$
\begin{array}{lll}
\Re(\underline{x})&=&(1+|\underline{x}|^2)(-1)^\ell2(\ell-\alpha)\underline{x} \omega_{\ell-\alpha-1}(\underline{x})\,\partial_{\underline{x}}^{\ell}\omega_{\alpha}(\underline{x})
\\&=&(-1)^\ell 2(\ell-\alpha)\underline{x} \omega_{\ell-\alpha}(\underline{x})\,\partial_{\underline{x}}^{\ell}\omega_{\alpha}(\underline{x}) \\
&=& \Im(\underline{x}).
\end{array}
$$
Hence,
$$
\begin{array}{lll}
G_{\ell+1,m}^{\alpha,\beta}(\underline{x})&=& -(1+|\underline{x}|^2)(-1)^{\ell}\omega_{\ell-\alpha}(\underline{x})\partial_{\underline{x}}^{(\ell+1)} \omega_{\alpha}(\underline{x})\\
&=& (-1)^{\ell+1}\omega_{\ell-\alpha+1}(\underline{x})\partial_{\underline{x}}^{(\ell+1)}\omega_{\alpha}(\underline{x}).
\end{array}
$$
\section{A 2-parameters Clifford-Jacobi polynomials and associated wavelets }
We propose in this section to introduce a 2-parameters class of polynomials based on Clifford-Jacobi ones. We denote such polynomials along the whole section by  $Z_{\ell,m}^{\alpha,\beta}(\underline{x})$. These are generated by the weight function
$$
\omega_{\alpha,\beta}(\underline{x})=(1-|\underline{x}|^2)^\alpha (1+|\underline{x}|^2)^\beta
$$
and its CK-extension $F^*$ expressed by
$$
\begin{array}{lll}
F^*(t,\underline{x})&=&
\displaystyle\sum\limits_{\ell=0}^{\infty}\dfrac{t^\ell}{\ell!}Z_{\ell,m}^{\alpha,\beta}(\underline{x})\,\omega_{\alpha-\ell,\beta-\ell}(\underline{x}).
\end{array}
$$
For more details about the definition of powers of the form $(1\pm\underline{u})^\alpha$ in Clifford analysis, we may refer to \cite{Brackx-Schepper-Sommen0}, \cite{Brackx-Schepper-Sommen1}, \cite{Brackx-Schepper-Sommen2}, \cite{Brackx-Schepper-Sommen3}, \cite{Brackx-Schepper-Sommen4} or \cite{Brackx-Schepper-Sommen5}. Next, we have in one hand
$$
\dfrac{\partial F^*(t,\underline{x})}{\partial  t}=\displaystyle\sum_{\ell=0}^{\infty}\dfrac{t^\ell}{\ell!} Z_{\ell+1,m}^{\alpha,\beta}(\underline{x})
\,\omega_{\alpha-\ell-1,\beta-\ell-1}(\underline{x}),
$$
and on the other hand,
$$
\begin{array}{lll}
\dfrac{\partial F^*(t,\underline{x})}{\partial \underline{x}}
&=&\displaystyle\sum_{\ell=0}^{\infty}\dfrac{t^\ell}{\ell!}\left( Z_{\ell,m}^{\alpha,\beta}(\underline{x}) \partial_{\underline{x}}\, \omega_{\alpha-\ell,\beta-\ell}(\underline{x})\right.\\
&&\qquad\qquad\quad\left.+\partial_{\underline{x}} (Z_{\ell,m}^{\alpha,\beta}(\underline{x}))\, \omega_{\alpha-\ell,\beta-\ell}(\underline{x}) \right),
\end{array}
$$
where
$$
\begin{array}{lll}
\partial_{\underline{x}} \omega_{\alpha-\ell,\beta-\ell}(\underline{x})
&=&-2(\alpha-\ell)\underline{x}\,\omega_{\alpha-\ell-1,\beta-\ell}(\underline{x})\\
&&\qquad\qquad+2(\beta-\ell)\underline{x}\,\omega_{\alpha-\ell,\beta-\ell-1}(\underline{x}).
\end{array}
$$
Then
$$
\begin{array}{lll}
\dfrac{\partial F^*(t,\underline{x})}{\partial \underline{x}}
&=&\displaystyle\sum_{\ell=0}^{\infty}\dfrac{t^\ell}{\ell!} Z_{\ell,m}^{\alpha,\beta}(\underline{x})
[-2(\alpha-\ell)\underline{x}\,\omega_{\alpha-\ell-1,\beta-\ell}(\underline{x})\\
&&\qquad+2(\beta-\ell)\underline{x}\,\omega_{\alpha-\ell,\beta-\ell-1}(\underline{x})]
+\partial_{\underline{x}} (Z_{\ell,m}^{\alpha,\beta}(\underline{x}))
\,\omega_{\alpha-\ell,\beta-\ell}(\underline{x}).
\end{array}
$$
From the monogenicity relation, we obtain
$$
\begin{array}{lll}
&&(\partial_t+\partial_{\underline{x}})F^*(t,\underline{x})\\
&=& Z_{\ell+1,m}^{\alpha,\beta}(\underline{x})
\,\omega_{\alpha-\ell-1,\beta-\ell-1}(\underline{x})
+Z_{\ell,m}^{\alpha,\beta}(\underline{x})
[-2(\alpha-\ell)\underline{x}\,\omega_{\alpha-\ell-1,\beta-\ell}\underline{x})\\
&&+2(\beta-\ell)\underline{x}\,\,\omega_{\alpha-\ell,\beta-\ell-1}(\underline{x})]
+\,\omega_{\alpha-\ell,\beta-\ell}(\underline{x})\,\partial_{\underline{x}}( Z_{\ell,m}^{\alpha,\beta}(\underline{x})) \\
&=& 0.
\end{array}
$$
Finally, we get the following result.
\begin{prop} The 2-parameters Clifford-Jacobi Polynomials $Z_{\ell,m}^{\alpha,\beta}$ satisfy the recurence relation
\begin{equation}\label{24}
\begin{array}{lll}
Z_{\ell+1,m}^{\alpha,\beta}(\underline{x})
&=&[2(\alpha-\ell)\underline{x}(1-\underline{x}^2)-2(\beta-\ell)\underline{x}\,(1+\underline{x}^2)] Z_{\ell,m}^{\alpha,\beta}(\underline{x})\\
&&-\,\omega_{1,1}(\underline{x})\partial_{\underline{x}} (Z_{\ell,m}^{\alpha,\beta}(\underline{x})).
\end{array}
\end{equation}
\end{prop}
\hskip-20pt For example, starting with $Z_{0,m}^{\alpha,\beta}(\underline{x}) =1$, a simple calculation yields that
$$
\begin{array}{lll}
Z_{1,m}^{\alpha,\beta}(\underline{x})
&=& 2\alpha\underline{x}(1-\underline{x}^2)-2\beta\underline{x} (1+\underline{x}^2)\\
&=&2(\alpha-\beta)\underline{x}-2(\alpha+\beta)\underline{x}^3
\end{array}
$$
For $\ell=1$, we get
$$
\begin{array}{lll}
&&Z_{2,m}^{\alpha,\beta}(\underline{x})\\
&=& [2(\alpha-1)\underline{x} (1-\underline{x}^2)-2(\beta-1)\underline{x}(1+\underline{x}^2)]
[2(\alpha-\beta)\underline{x}-2(\alpha+\beta)\underline{x}^3]
\\
&&+(1-\underline{x}^4) [-2(\alpha-\beta)m+2(\alpha+\beta)(m+2)\underline{x}^2]
\\
&=& 2(\alpha-\beta)m+[4\alpha(\alpha-1)+4\beta(\beta-1)-8\alpha\beta-2(\alpha+\beta)m]\underline{x}^2\\
&&+[8\beta(\beta-1)-8\alpha(\alpha-1)+2(\beta-\alpha)m]\underline{x}^4\\
&&+[4\alpha(\alpha-1)+4\beta(\beta-1)+8\alpha\beta+2(\beta+\alpha)m]\underline{x}^6.
\end{array}
$$
For $\ell=2$, we obtain
$$
\begin{array}{lll}
&&Z_{3,m}^{\alpha,\beta}(\underline{x})\\
&=&[[4(\alpha-\beta)^2-4(\alpha+\beta)]m+8\alpha(\alpha-1)+8\beta(\beta-1)-8\alpha\beta]\underline{x}\\
&+& [-26\alpha(\alpha-1)+40\beta(\beta-1)-16\alpha\beta-4(\alpha+\beta)+4(\beta-\alpha)[\alpha+\beta-2]]\underline{x}^3\\
&+& [16\beta(\beta-1)(\alpha-\beta)-16\alpha(\alpha-1)(\alpha-\beta)-4(\alpha-\beta)^2m]\\
&-&[4\alpha(\alpha-1)+4\beta(\beta-1)-8\alpha\beta-2(\alpha+\beta)m][2(\alpha+\beta-2)] \underline{x}^5\\
&+& [8\alpha(\alpha-1)+8\beta(\beta-1)+16\alpha\beta+4(\beta+\alpha)m](\alpha-\beta)\\
&-&2[8\beta(\beta-1)-8\alpha(\alpha-1)+2(\beta-\alpha)m](\alpha+\beta)\underline{x}^7\\
&-& [4\alpha(\alpha-1)+4\beta(\beta-1)+8\alpha\beta+2
(\beta+\alpha)m[2\alpha+2\beta-2] \underline{x}^9.
\end{array}
$$
Remark that $Z_{\ell,m}^{\alpha,\beta}(\underline{x})$ is a polynomial of degree $3\ell$ in $\underline{x}$.
\begin{prop} The 2-parameters Clifford-Jacobi polynomials $Z_{\ell,m}^{\alpha,\beta}$ may be obtained via the Rodrigues formula
\begin{equation}\label{25}
Z_{\ell,m}^{\alpha,\beta}(\underline{x})=(-1)^\ell \,\omega_{\ell-\alpha,-\ell-\beta}(\underline{x}) \,
\partial_{\underline{x}}^{\ell} [(1+\underline{x}^2)^\alpha (1-\underline{x}^2)^\beta].
\end{equation}
\end{prop}
\hskip-20pt\textbf{Proof.} For $\ell=0$, the situation is obvious. For $\ell=1$, we have
$$
\begin{array}{lll}
\partial_{\underline{x}} ( \,\omega_{\alpha,\beta}(\underline{x}))
&=&-2\alpha\underline{x}\,\omega_{\alpha-1,\beta}(\underline{x})+2\beta\underline{x}\,\omega_{\alpha,\beta-1}(\underline{x})\\
&=&(-1)\,\,\omega_{\alpha-1,\beta-1}(\underline{x})[2\alpha\underline{x}(1-\underline{x}^2)-2\beta\underline{x}(1+\underline{x}^2)]\\
&=&(-1)\,\,\omega_{\alpha-1,\beta-1}(\underline{x})Z_{1,m}^{\mu,\alpha}(\underline{x}).
\end{array}
$$
Thus,
$$
Z_{1,m}^{\alpha,\beta}(\underline{x})=(-1)\omega_{1-\alpha,1-\beta}(\underline{x})\partial_{\underline{x}}(\omega_{\alpha,\beta}(\underline{x})).
$$
For $\ell=2$, we get
$$
\begin{array}{lll}
&&
\partial_{\underline{x}}^2(\omega_{\alpha,\beta}(\underline{x}))\\
&=&-2\alpha[\partial_{\underline{x}}(\underline{x}(1+\underline{x}^2)^{\alpha-1})(1-\underline{x}^2)^\beta
+\underline{x}(1+\underline{x}^2)^{\alpha-1}\,\partial_{\underline{x}}(1-\underline{x}^2)^{\beta}]\\
&&+ 2\beta[\partial_{\underline{x}}(\underline{x}(1+\underline{x}^2)^{\alpha})(1-\underline{x}^2)^{\beta-1}
+\underline{x}(1+\underline{x}^2)^{\alpha}\partial_{\underline{x}}(1-\underline{x}^2)^{\beta-1}]\\
&=& 2m\alpha\omega_{\alpha-1,\beta}(\underline{x})+4\alpha(\alpha-1)\underline{x}^2\omega_{\alpha-2,\beta}(\underline{x})\\
&&- 4\alpha\beta\underline{x}^2\omega_{\alpha-1,\beta-1}(\underline{x})-2m\beta\,\omega_{\alpha,\beta-1}(\underline{x})\\
&&- 4\alpha\beta\underline{x}^2\omega_{\alpha-1,\beta-1}(\underline{x})+4\beta(\beta-1)\underline{x}^2\omega_{\alpha,\beta-2}(\underline{x})\\
&=&(-1)^2\omega_{\alpha-2,\beta-2}(\underline{x})[2m\alpha\omega_{1,2}(\underline{x})\\
&&+4\alpha(\alpha-1)\underline{x}^2(1-\underline{x}^2)^2-4\alpha\beta\underline{x}^2\omega_{1,1}(\underline{x})\\
&&-2m\beta\omega_{2,1}(\underline{x})-4\alpha\beta\underline{x}^2\omega_{1,1}(\underline{x})+4\beta(\beta-1)\underline{x}^2(1+\underline{x}^2)^2]\\
&=&(-1)^2\omega_{\alpha-2,\beta-2}(\underline{x})Z_{2,m}^{\alpha,\beta}(\underline{x}).
\end{array}
$$
Then
$$
Z_{2,m}^{\alpha,\beta}(\underline{x})=(-1)^2\omega_{2-\alpha,2-\beta}(\underline{x})\partial_{\underline{x}}(\omega_{\alpha,\beta}(\underline{x}).
$$
Now assume that
$$
Z_{\ell,m}^{\alpha,\beta}(\underline{x})=(-1)^\ell \,\omega_{\ell-\alpha,\ell-\beta}(\underline{x}) \,
\partial_{\underline{x}}^{\ell} \,\omega_{\alpha,\beta}(\underline{x}).
$$
Denote
$$
\begin{array}{lll}
\wp&=&[2(\alpha-\ell)\underline{x}(1-\underline{x}^2)-2(\beta-\ell)\underline{x}\,(1+\underline{x}^2)]\\
&&\qquad\qquad\qquad[(-1)^\ell \,\omega_{\ell-\alpha,\ell-\beta}(\underline{x}) \,
\partial_{\underline{x}}^{\ell} \,\omega_{\alpha,\beta}(\underline{x})]
\end{array}
$$
and
$$
\begin{array}{lll}
\aleph&=&(-1)^\ell\,\omega_{1,1}(\underline{x})\,
[2(\alpha-\ell)\underline{x}\,\omega_{\ell-\alpha-1,\ell-\beta}(\underline{x})\\
&&\qquad\qquad\qquad-2(\beta-\ell)\underline{x}\,\,\omega_{\ell-\alpha,\ell-\beta-1}(\underline{x})]	\partial_{\underline{x}}^{\ell} \,\omega_{\alpha,\beta}(\underline{x}) .
\end{array}
$$
Then we derive from (\ref{24}) and (\ref{25}) that
$$
\begin{array}{lll}
\;&\;& Z_{\ell+1,m}^{\alpha,\beta}(\underline{x})\\
&=& [2(\alpha-\ell)\underline{x}(1-\underline{x}^2)-2(\beta-\ell)\underline{x}\,(1+\underline{x}^2)]Z_{\ell,m}^{\alpha,\beta}(\underline{x})\\
&&\qquad\qquad\qquad-\,\omega_{1,1}(\underline{x})\partial_{\underline{x}} Z_{\ell,m}^{\alpha,\beta}(\underline{x}) \\
&=&\wp-\aleph-(-1)^\ell \,\omega_{\ell-\alpha+1,\ell-\beta+1}(\underline{x})\,\partial_{\underline{x}}^{\ell+1}\omega_{\alpha,\beta}(\underline{x}).
\end{array}
$$
Otherwise, we have
$$
\begin{array}{lll}
\aleph&=& (-1)^\ell [2(\alpha-\ell)\underline{x}\,\omega_{\ell-\alpha,\ell-\beta+1}(\underline{x})\\
&&\qquad\qquad\qquad-2(\beta-\ell)\underline{x}\,\,\omega_{\ell-\alpha+1,\ell-\beta-1}(\underline{x})]
\partial_{\underline{x}}^{\ell}\omega_{\alpha,\beta}(\underline{x})\\
&=& (-1)^\ell[2(\alpha-\ell)\underline{x}(1-\underline{x}^2)\\
&&\qquad\qquad\qquad-2(\beta-\ell)\underline{x}(1+\underline{x}^2)]\,\omega_{\ell-\alpha,\ell-\beta}(\underline{x}) \partial_{\underline{x}}^{\ell}\omega_{\alpha,\beta}(\underline{x})\\
&=& \wp.
\end{array}
$$
Hence,
$$
Z_{\ell+1,m}^{\alpha,\beta}(\underline{x})= (-1)^{\ell+1} \,\omega_{\ell-\alpha+1,\ell-\beta+1}(\underline{x})\,\partial_{\underline{x}}^{\ell+1}(\omega_{\alpha,\beta}(\underline{x})).
$$
The following orthogonality relation is proved.
\begin{prop}\label{ZlmOrthogonality}
Let
$$
I_{\ell,t,p}^{\alpha,\beta}=\displaystyle\int_{\mathbb{R}^m}\underline{x}^{\ell} Z_{t,m}^{\alpha+p,\beta+p}(\underline{x})\, \omega_{\alpha,\beta}(\underline{x})\, dV(\underline{x}).
$$
For $4t<1-m-2(\alpha+\beta)$ we have
\begin{equation}\label{26}
I_{\ell,t,t}^{\alpha,\beta}=0 .
\end{equation}	
\end{prop}
\hskip-20pt\textbf{Proof.} Denote
$$
I_{\ell,t}= \displaystyle\int_{\mathbb{R}^m} \underline{x}^{\ell}\,\partial_{\underline{x}}^t( \omega_{\alpha+t,\beta+t}(\underline{x})\,dV(\underline{x})).
$$
Using Stokes's theorem, we obtain
$$
\begin{array}{lll}
&&\displaystyle\int_{\mathbb{R}^m}\underline{x}^{\ell} Z_{t,m}^{\alpha+t,\beta+t}(\underline{x})\omega_{\alpha,\beta}(\underline{x})dV(\underline{x})\\
&=&\displaystyle\int_{\mathbb{R}^m} \underline{x}^{\ell} (-1)^t \omega_{t-\alpha-t,t-\beta-t}(\underline{x}) \partial_{\underline{x}}^t(\omega_{\alpha+t,\beta+t}(\underline{x}))\omega_{\alpha,\beta}(\underline{x}) \,dV(\underline{x})\\
&=& (-1)^t\displaystyle\int_{\mathbb{R}^m} \underline{x}^{\ell}\,\partial_{\underline{x}}^t( \omega_{\alpha+t,\beta+t}(\underline{x})\,dV(\underline{x}))\\
&=& (-1)^t\displaystyle\int_{\mathbb{R}^m} \underline{x}^{\ell}\,\partial_{\underline{x}}\partial_{\underline{x}}^{t-1}(\omega_{\alpha+t,\beta+t}(\underline{x}))\,dV(\underline{x}))\\
&=&(-1)^t\left[\displaystyle\int_{\partial\mathbb{R}^m}\underline{x}^{\ell}\partial_{\underline{x}}^{t-1}\omega_{\alpha+t,\beta+t}(\underline{x})\,dV(\underline{x})\right.\\
&&\qquad\qquad\qquad
\left.-\displaystyle\int_{\mathbb{R}^m}\partial_{\underline{x}}(\underline{x}^{\ell})\partial_{\underline{x}}^{t-1}\omega_{\alpha+t,\beta+t}(\underline{x})dV(\underline{x})\right].
\end{array}
$$
Denote already
$$
I=\displaystyle\int_{\partial\mathbb{R}^m}\underline{x}^{\ell}\partial_{\underline{x}}^{t-1}\omega_{\alpha+t,\beta+t}(\underline{x})\,dV(\underline{x})
$$
and
$$
II=\displaystyle\int_{\mathbb{R}^m}\partial_{\underline{x}}(\underline{x}^{\ell})\partial_{\underline{x}}^{t-1}\omega_{\alpha+t,\beta+t}(\underline{x})dV(\underline{x}).
$$
The integral $I$ vanishes due to the assumption
$$
0<t<\dfrac{1-m-2(\alpha+\beta)}{4}.
$$
Due to Lemma \ref{lemmaevenodd}, the second satisfies
$$
II= \gamma_{l,m}\displaystyle\int_{\mathbb{R}^m}\underline{x}^{\ell-1}
\partial_{\underline{x}}^{t-1}\omega_{\alpha+t,\beta+t}(\underline{x})\,dV(\underline{x})=\gamma_{l,m} I_{\ell-1,t-1}.
$$
Hence, we obtain
$$
\begin{array}{lll}
&&\displaystyle\int_{\mathbb{R}^m}\underline{x}^\ell Z_{t,m}^{\alpha+t,\beta+t}(\underline{x})\omega_{\alpha,\beta}(\underline{x}) dV(\underline{x})\\
&=&(-1)^{t+1} \gamma_{l,m} I_{\ell-1,t-1}\\
&=&  (-1)^{t+1}  \gamma_{l,m}[ (-1)^{t} \gamma_{l-1,m}I_{\ell-2,t-2}]\\
&=& (-1)^{2t+1} \gamma_{l,m}  \gamma_{l-1,m} I_{\ell-2,t-2}\\
&\vdots&\\
&=&C(m,\ell,t) I_0\\
&=& 0.
\end{array}
$$
where $C(m,\ell,t)=(-1)^{ml+1}\displaystyle\prod_{k=0}^{m}\gamma_{k,m}$.
\begin{defn} The generalized 2-parameters Clifford-Jacobi wavelet mother is defined by
$$
\psi_{\ell,m}^{\alpha,\beta}(\underline{x})
=Z_{\ell,m}^{\alpha+\ell,\beta+\ell}(\underline{x})
\omega_{\alpha,\beta}(\underline{x})
=(-1)^\ell\partial_{\underline{x}}^{(\ell)}\omega_{\alpha+\ell,\beta+\ell}(\underline{x}).
$$
\end{defn}
\hskip-20pt Furthermore, the wavelet $\psi_{\ell,m}^{\alpha,\beta}(\underline{x})$ have vanishing moments as is shown in the next proposition.
\begin{prop} The following assertions hold.
\begin{enumerate}
\item For $0<k<-m-\ell-2(\alpha+\beta) $ and $k<\ell$ we have
\begin{equation}\label{27}
\displaystyle\int_{\mathbb{R}^m} \underline{x}^k \psi_{\ell,m}^{\alpha,\beta}(\underline{x}) dV(\underline{x})=0.
\end{equation}
\item Its Clifford-Fourier transform is
\begin{equation}
\widehat{\psi_{\ell,m}^{\alpha,\beta}(\underline{u})}=(-i)^\ell\,\underline{\xi}^\ell(2\pi)^{\frac{m}{2}}\rho^{1-\frac{m}{2}+\ell}\,\displaystyle\int_{0}^\infty\widetilde{\omega}_{\alpha,\beta}^l(r)\,J_{\frac{m}{2}-1}(r\rho)dr.
\end{equation}
where
$$
\widetilde{\omega}_{\alpha,\beta}^l(r)=((1-r^2)\varepsilon_r)^{\alpha+\ell} (1+r^2)^{\beta+\ell} r^{\frac{m}{2}}
$$
with $\varepsilon_r=\mbox{sign}(1-r)$.
\end{enumerate}
\end{prop}
\hskip-20pt\textbf{Proof.} The first assertion is a natural consequence of Proposition \ref{ZlmOrthogonality}. We prove the second. We have
$$
\begin{array}{lll}
\widehat{\psi}_{\ell,m}^{\alpha,\beta}(\underline{u})
&=&\displaystyle\int_{\mathbb{R}^m}\psi_{\ell,m}^{\alpha,\beta}(\underline{x})e^{-i\underline{x}.\underline{u}}\,dV(\underline{x})\\
&=&(-1)^\ell\displaystyle\int_{\mathbb{R}^m}\partial_{\underline{x}}^\ell\left(\omega_{\alpha+\ell,\beta+\ell}(\underline{x})\right)e^{-i\underline{x}.\underline{u}}\,dV(\underline{x})\\
&=&(-1)^\ell\displaystyle\int_{\mathbb{R}^m}\omega_{\alpha+\ell,\beta+\ell}(\underline{x})e^{-i\underline{x}.\underline{u}}(i\underline{u})^\ell \,dV(\underline{x})\\
&=&(-1)^\ell\,(i\underline{u})^\ell\displaystyle\int_{\mathbb{R}^m}\omega_{\alpha+\ell,\beta+\ell}(\underline{x})\, e^{-i\underline{x}.\underline{u}}\,dV(\underline{x})\\
&=&(-1)^{\ell}\,(i\underline{u})^{\ell}\displaystyle\int_{\mathbb{R}^m}
(1-|\underline{x}|^2)^{\alpha+\ell}(1+|\underline{x}|^2)^{\beta+\ell}\,e^{-i\underline{x}.\underline{u}} \,dV(\underline{x})\\
&=&(-1)^\ell(i\underline{u})^\ell\widehat{\omega_{\alpha+\ell,\beta+\ell}}(\underline{u}).\end{array}
$$
This Fourier transform  can be simplified by using the spherical co-ordinates. By definition, we have
\begin{equation}\label{spherical-co-ordinates}
\widehat{\omega_{\alpha+\ell,\beta+\ell}}(\underline{u})=\displaystyle\int_{\mathbb{R}^m} (1-|\underline{x}|^2)^{\alpha+\ell} (1+|\underline{x}|^2)^{\beta+\ell}\, e^{-i<\underline{x},\underline{u}>} dV(\underline{x})
\end{equation}
Introducing spherical co-ordinates
$$
\underline{x}=r\underline{\omega},\quad \underline{u}=\rho\underline{\xi},\quad r=|\underline{x}|,\quad \rho=|\underline{u}|, \quad \underline{\omega}\in S^{m-1}, \,\underline{\xi}\in S^{m-1}
$$
expression (\ref{spherical-co-ordinates}) becomes
$$
\begin{array}{lll}
\widehat{\omega_{\alpha+\ell,\beta+\ell}(\underline{u})}&=&\displaystyle\int_{0}^\infty\widetilde{\omega}_{\alpha,\beta}^l(r)\,r^{\frac{m}{2}-1}\,dr\displaystyle\int_{S^{m-1}} e^{-i<r\underline{\omega},\rho\underline{\xi}>} d\sigma(\underline{\omega})
\end{array}
$$
where $d\sigma(\underline{\omega})$ stands for the Lebesgue measure on $S^{m-1}$.\\
We now use the following technical result which is known in the theory of Fourier analysis of radial functions and the theory of bessel functions.
\begin{lem}\label{BesselFourierTransform}
$$
\displaystyle\int_{S^{m-1}}e^{-i<r\underline{\omega},\rho\underline{\xi}>} d\sigma(\underline{\omega})=\displaystyle\frac{(2\pi)^{\frac{m}{2}} J_{\frac{m}{2}-1} (r\rho)}{r^{\frac{m}{2}-1} \rho^{\frac{m}{2}-1}}
$$
where $J_{\frac{m}{2}-1}$ is the bessel function of the first kind of order $\frac{m}{2}-1$.
\end{lem}
The proof of this result is a good exercice in Fourier analysis and orthogonal transformations and may be found in \cite{Stein-Weiss}.\\
Now, according to Lemma \ref{BesselFourierTransform}, we obtain
$$
\begin{array}{lll}
\widehat{\omega_{\alpha+\ell,\beta+\ell}(\underline{u})}&=&(2\pi)^{\frac{m}{2}}\rho^{1-\frac{m}{2}+\ell}\,\displaystyle\int_{0}^\infty \widetilde{\omega}_{\alpha,\beta}^l(r)J_{\frac{m}{2}-1} (r\rho) dr.
\end{array}
$$
Consequently, we obtain the following expression for the Fourier transform of the $(\alpha,\beta)$-Clifford-Jacobi wavelets
$$
\begin{array}{lll}
\widehat{\psi_{\ell,m}^{\alpha, \beta}(\underline{u})}&=&(-i)^\ell\,\xi^\ell (2\pi)^{\frac{m}{2}}\rho^{1-\frac{m}{2}+\ell}\,\displaystyle\int_{0}^\infty \widetilde{\omega}_{\alpha,\beta}^l(r)J_{\frac{m}{2}-1} (r\rho) dr.
\end{array}
$$
\begin{defn}
The copy of the generalized 2-parameters Clifford-Jacobi wavelet at the scale $a>0$ and the position $\underline{b}$ is defined by
$$
_a^{\underline{b}}\psi_{\ell,m}^{\alpha,\beta}(\underline{x})=a^{-\frac{m}{2}}\psi_{\ell,m}^{\alpha,\beta}(\dfrac{\underline{x}-\underline{b}}{a}).
$$
\end{defn}
\begin{defn}
The wavelet transform of a function $f$ in $L_2$ according to th generalized 2-parameters Clifford-Jacobi wavelet at the scale $a$ and the position $\underline{b}$ is
$$
C_{a,\underline{b}}(f)=<f,\,_a^{\underline{b}}\psi_{\ell,m}^{\alpha,\beta}>=
\displaystyle\int_{\mathbb{R}^m}f(x)\,_a^{\underline{b}}\psi_{\ell,m}^{\alpha,\beta}(\underline{x})dV(\underline{x}).
$$
\end{defn}
\hskip-20pt The following Lemma guarantees that the candidate $\psi_{\ell,m}^{\alpha,\beta}$ is indeed a mother wavelet. Analogue result is already checked in \cite{Brackx-Schepper-Sommen1}. The proof is based on the asymptotic behaviour of Bessel functions and thus left to the reader.
\begin{lem}\label{Admissibilityofthenewwavelets}
The quantity
$$
\mathcal{A}_{\ell,m}^{\alpha,\beta}=\dfrac{1}{\omega_m}\displaystyle\int_{\mathbb{R}^m}|\widehat{\psi_{\ell,m}^{\alpha,\beta}}(\underline{x})|^2 \dfrac{dV(\underline{x})}{|\underline{x}|^m}
$$
is finite. ($\omega_m$ is the volume of the unit sphere $S^{m-1}$ in $\mathbb{R}^m$.
\end{lem}
To state the final result dealing with the reconstruction formula relatively to the constructed new wavelets, we introduce firstly the inner product
$$
<C_{a,\underline{b}}(f),C_{a,\underline{b}}(g)>=\dfrac{1}{\mathcal{A}_{\ell,m}^{\alpha,\beta}}\displaystyle\int_{\mathbb{R}^m}\displaystyle\int_{0}^{+\infty}\overline{C_{a,\underline{b}}(f)}C_{a,\underline{b}}(g)\dfrac{da}{a^{m+1}}dV(\underline{b}).
$$
We obtain the following result.
\begin{thm}\label{ReconstructionFormula}
Then, the function $f$ may be reconstructed by
$$
f(x)=\dfrac{1}{\mathcal{A}_{\ell,m}^{\alpha,\beta}}\displaystyle\int_{a>0}\displaystyle\int_{b\in\mathbb{R}^m}C_{a,\underline{b}}(f)\,\psi\left(\dfrac{\underline{x}-\underline{b}}{a}\right)\dfrac{da\,dV(\underline{b})}{a^{m+1}}.
$$
where the equality has to be understood in the $L_2$-sense.
\end{thm}
\hskip-20pt The proof reposes on the following result.
\begin{lem}\label{ProduitScalaireCoefficient}
It holds that
$$
\displaystyle\int_{a>0}\displaystyle\int_{b\in\mathbb{R}^m}\overline{C_{a,\underline{b}}(f)}C_{a,\underline{b}}(g)\displaystyle\frac{da\,dV(\underline{b})}{a^{m+1}}
=\mathcal{A}_{\ell,m}^{\mu,\alpha}\displaystyle\int_{\mathbb{R}^m}f(\underline{x})\overline{g(\underline{x})}dV(\underline{x}).
$$
\end{lem}
\hskip-20pt\textbf{Proof.} Using the Clifford Fourier transform we observe that
$$
C_{a,\underline{b}}(f)(\underline{b})=\widetilde{a^{\frac{m}{2}}\widehat{\widehat{f}(\underline{.})\widehat{\psi}(a\underline{.})}}(\underline{b}),
$$
where, $\widetilde{h}(\underline{u})=h(-\underline{u})$, $\forall\,h$. Thus,
$$
\overline{C_{a,\underline{b}}(f)}C_{a,\underline{b}}(g)=\overline{\widehat{\left(\widehat{f}(\underline{.})a^{\frac{m}{2}} \widehat{\psi}(a\underline{.})\right)}}(-\underline{b}) \widehat{\left(\widehat{g}(\underline{.})a^{\frac{m}{2}} \widehat{\psi}(a\underline{.})\right)}(-\underline{b}).
$$
Consequently,
$$
\begin{array}{lll}
<C_{a,\underline{b}}(f),C_{a,\underline{b}}(g)>
&=&\displaystyle\int_{a>0}\displaystyle\int_{\mathbb{R}^m}\overline{\widehat{\widehat{f}(\underline{.})a^{\frac{m}{2}}\widehat{\psi}(a\underline{.})}}\,\widehat{\widehat{g}(\underline{.})a^{\frac{m}{2}}\widehat{\psi}(a\underline{.})}\displaystyle\frac{da\,dV(\underline{b})}{a^{m+1}}\\
&=& \displaystyle\int_{a>0}\displaystyle\int_{\mathbb{R}^m} \overline{\widehat{f}(\underline{b})}\widehat{g}(\underline{b})
\displaystyle\frac{a^m|\widehat{\psi}(a\underline{b})|^2}{a^{m+1}} \,da\,dV(\underline{b})\\
&=&{\mathcal{A}_{\ell,m}^{\mu,\alpha}}\displaystyle\int_{\omega\in\mathbb{R}^m}\overline{\widehat{f(\underline{b})}}\widehat{g}(b)dV(\underline{b})\\
&=&{\mathcal{A}_{\ell,m}^{\mu,\alpha}}<\widehat{f},\widehat{g}>\\
&=& <f,g> .
\end{array}
$$
\textbf{Proof of Theorem \ref{ReconstructionFormula}.} It follows immediately from lemma \ref{ProduitScalaireCoefficient}.
\section{Back to Legendre and Chebyshev polynomials}
This section is devoted to multi folds. One of the aims is to discuss the link of the present work and the possibility to construct Clifford-Legender and Clifford-Techebyshev wavelets. Such an aim is itself a motivation to show the role of the parameters appearing in the present construction and the fact that it consists of a large class of polynomials and consequently wavelets that may englobe for instence a Legendre and a Techebyshev case. Befor doing that we stress on the fact that in our knwledge there are no previous works that have developed the special cases of Legendre and Techebyshev polynomials in the context of Clifford analysis. On work found is dealing with general orthogonal polynomials in the Clifford context reminiscent of some brief comming back to the explicit expression of Legendre polynomials in the real case.

Recall firstly that when relaxing the parameter $\alpha$ for the weight function applied in the present work, we get the Gegenbauer existing case developed in \cite{Brackx-Schepper-Sommen1} and \cite{DeSchepper}. Recall also that classical real analysis affirm that Legendre, Techebyshev and Gegenbauer polynomials are issued from three-level recurence relations such as
\begin{equation}\label{LegendrePolynomials}
L_{n+1}=\dfrac{2n+1}{n+1}XL_{n}-\dfrac{n}{n+1}L_{n-1},\quad \forall n \in\mathbb{N}^{*}
\end{equation}
for Legendre polynomials,
\begin{equation}\label{TchebychevPolynomials}
T_{n+1}=2XT_{n}-T_{n-1},\quad \forall n \in\mathbb{N}^{*}
\end{equation}
for Tchebyshev polynomials and
\begin{equation}\label{GegenbauerPolynomials}
mG_{m}^{p}(x)=2x(m+p-1)G_{m-1}^{p}(x)-(m+2p-2)G_{m-2}^{p}(x),
\end{equation}
for Gegenbauer ones. These three classes of polynomials are defined also by means of the Rodrigues rule. Legendre polynomials are given by
\begin{equation}\label{LegendrePolynomialsRodrigues}
L_{n}(x)=\dfrac{d^n}{dx^{n}}[\dfrac{(x^{2}-1)^{n}}{2^{n}n!}]
\end{equation}
which with the Leibnitz rule yields an explicit form
\begin{equation}\label{LegendrePolynomialsExplicit}
L_{n}(x)=\dfrac{1}{2^{n}}\displaystyle\sum_{k=0}^n(C_n^k)^2(x-1)^{n-k}(x+1)^k.
\end{equation}
Tchebyshev polynomials are expressed via the Rodrigues rule as
\begin{equation}\label{ChebyshevPolynomialsRodrigues}
T_{n}(x)=\dfrac{2^n(-1)^nn!}{(2n)!}(1-x^{2})^{\frac{1}{2}}\dfrac{d^{n}}{dx^{n}}((1-x^{2})^{n-\frac{1}{2}}),
\end{equation}
which also with the Leibnitz rule induces the explicit form
\begin{equation}\label{TchebyshevPolynomialsExplicit}
T_{n}(x)=\dfrac{1}{2^{n}}\displaystyle\sum_{k=0}^nC_{2n}^{2k}(x-1)^{n-k}(x+1)^k.
\end{equation}
Gegenbauer polynomials called also ultra-spheroidal polynomials may be introduced via the Rodrigues rule as
\begin{equation}\label{GegenbauerPolynomialsRodrigues}
G^{p}_{m}(x)=(-1)^{m}\omega_{m,p}(1-x^{2})^{\frac{1}{2}-p}\dfrac{d^{m}}{dx^{m}}\left((1-x^{2})^{p+m-\frac{1}{2}}\right),
\end{equation}
where
$$
\omega_{m,p}=\dfrac{2^{m}m!\Gamma(p+\frac{1}{2})\Gamma(m+2p)}{\Gamma(2p)\Gamma(p+m+\frac{1}{2})}.
$$
Applying already the Leibnitz derivation rule, we obtain
\begin{equation}\label{GegenbauerPolynomialsExplicit}
G^{p}_{m}(x)=\dfrac{\omega_{m,p}}{2^{2m}m!}\displaystyle\sum_{k=0}^mC_{2m}^{2k}(x-1)^{m-k}(x+1)^k.
\end{equation}
Relatively to weight functions, Legendre polynomials are related to the weight $\omega_L(x)=(1-x^2)^n$. Tchebyshev polynomials are issued from the weight function $\omega_T(x)=(1-x^2)^{-1/2}$. Gegenbauer polynomials are deduced from the weight function $\omega_G(x)=(1-x^2)^{p-\frac{1}{2}}$. This means that both Legendre plynomials and Tchebyshev ones may be deduced from Gegenbauer case by simple choices of the parameter $p$ in $\omega_G$, such that $p=n+\frac{1}{2}$ for Legendre polynomials and $p=0$ for Techebyshev ones. So, a motivated extension in the case of Clifford analysis may be conducted by applying the same values for the parameter $p$ in Clifford-Gegenbauer polynomials explicit form to obtain explicit forms for Clifford-Legendre and Clifford-Techebyshev extensions. Therefore, an eventual form for Clifford-Legendre polynomials will be deduced from (\ref{LegendrePolynomialsExplicit}) as
\begin{equation}\label{CliffordLegendrePolynomialsExplicit}
L_{n}(\underline{x})=\dfrac{1}{2^{n}}\displaystyle\sum_{k=0}^n(C_n^k)^2(-1)^{n-k}(1-\underline{x})^{n-k}(1+\underline{x})^k
\end{equation}
and similarly an explicit Clifford extension of Techebyshev polynomials may be
\begin{equation}\label{CliffordTchebyshevPolynomialsExplicit}
T_{n}(\underline{x})=\dfrac{1}{2^{n}}\displaystyle\sum_{k=0}^nC_{2n}^{2k}(-1)^{n-k}(1-\underline{x})^{n-k}(1+\underline{x})^k.
\end{equation}
A second idea may be developed by comparing the induction rules in the real case and adopt similar changes on the Gegenbauer induction rule in Clifford analysis to obtain induction rules for Legendre and Tchebyshev ones. The explicit forms (\ref{CliffordLegendrePolynomialsExplicit}) and (\ref{CliffordTchebyshevPolynomialsExplicit}) may be good starting points to explicit the induction rules expected. We acheive our paper by proposing these facts as future works.

As a result of these facts, it is immediate that the present case englobe a large set of Clifford polynomials and consequently of Clifford wavelets that may be adopted to the cited cases above.
\section{Conclusion}
In this paper we have introduced new classes of orthogonal polynomials relatively to new 2-parameters weight in the context of Clifford analysis. The new class generalizes the well known Jacobi and Gegenbauer polynomials. Such polynomial sets are next applied to introduce new wavelets in Clifford analysis. Fourier-Plancherel type results are proved for the new classes.

\end{document}